\documentclass[11pt,twoside]{article}
\usepackage{amsmath,amsbsy,amsfonts,amsthm}
\newtheorem{thm}{Theorem}[section]
\newtheorem{prop}{Proposition}[section]
\newtheorem{lem}{Lemma}[section]
\newtheorem{rmk}{Remark}[section]

\newtheorem{cor}{Corollary}[section]
\newtheorem{clm}{Claim}[section]

\newcommand{\eproof}{\hfill \mbox{${\square}$}}

\def\proof{\noindent {\bf Proof.}}

\def\RR{\mathbb{R}}
\def\NN{\mathbb{N}}

\def \o{o_n(1)}

\def \ep {\epsilon}
\def \la {\lambda}
\def \e {\epsilon}

\theoremstyle{definition}
\begin{document}

\title{Multiplicity results for a class of quasilinear equations with exponential critical growth}
\author{Claudianor O. Alves and Luciana R. de Freitas}

\date{}
\maketitle

\begin{abstract}
In this work, we prove the existence and multiplicity of positive solutions for the following class of quasilinear elliptic equations

$$
\left \{
\begin{array}{lll}
-\epsilon^{N}\Delta_{N} u + \left(1+\mu A(x) \right)\left| u\right|^{N-2}u= f(u)\,\,\,\, \mbox{ in }  \,  \RR^{N},\\
u > 0 \,\,\,\, \mbox{ in }  \, \RR^{N},
\end{array}
\right.
$$
 \\
where $\Delta_{N}$ is the N-Laplacian operator, $N \geq 2$,  $f$ is a function with exponential critical growth, $\mu$ and $\epsilon$ are positive parameters and $A$ is a nonnegative continuous function verifying some hypotheses. To obtain our results, we combine variational arguments and  Lusternik-Schnirelman category theory .

\end{abstract}

\vspace{0.5 cm}



\section{Introduction}
\hspace{0.6 cm}In this article, we consider the following class of quasilinear problem
$$
\left \{
\begin{array}{lll}
-\epsilon^N\Delta_{N} u + (1+\mu A(x))|u|^{N-2}u = f(u)\,\,\,\, \mbox{ in }  \, \RR^N (N\geq2),\\
\quad u > 0\,\,\,\, \mbox{ in }  \, \RR^{N}, \\
\end{array}
\right.\eqno{ (P_{\mu,\epsilon})}
$$
where $\Delta_{N}u= div( \left|\nabla u\right|^{N-2}\nabla u ) $ is the $N$-Laplacian operator, $\mu$ and $\epsilon$ are positive parameters and the nonlinear term $f$ is a function having critical exponential growth. The hypotheses on $A$ are the same assumed in \cite{CS} , namely:
\begin{description}
		\item[${(A_1)}$] $A \in C^1(\RR^N,\RR)$ is a nonnegative function such that $A^{-1}(0)=\overline{\Omega}\cup D$, where $\Omega = \,$ $ intA^{-1} (0)$  is a bounded open set with smooth boundary, $0 \in \Omega$ and $m (D)=0$;  
		\item[${(A_2)}$ ] There exists $M_0 >0$ such that $m(\{x \in \RR^N ; A(x) \leq M_0\}) < \infty$,
\end{description}
where $m$ denotes the Lebesgue measure on $\RR^{N}$. From $(A_1)$, we can fix $r>0$ such that $B_r(0) \subset \Omega$ and the sets
\[
\Omega^+ = \left\{x \in \mathbb{R}^N \, ; \, d(x, \overline{\Omega}) \leq r\right\}
\]
and
\[
\Omega^- = \left\{x \in \Omega \, ; \, d(x, \partial \Omega) \geq r\right\}
\]
are homotopically equivalent to $\Omega$.

The hypotheses on the nonlinear term $f$ are as follows:  
\begin{description}
    \item[$(H_0)$] $f \in C^{1}(\mathbb{R},\mathbb{R})$ is a function with exponential critical growth, that is, there exists  $\alpha_{0} > 0$ such that
    $$
    \lim_{\left|s\right| \rightarrow \infty} \frac{\left|f(s)\right|}{e^{\alpha\left|s\right|^{\frac{N}{N-1}}}}=
     \left\{
\begin{array}{lll}
0, & \mbox{ if } &\alpha > \alpha_{0}, \\
+\infty, &  \mbox{ if } & \alpha < \alpha_{0}.
\end{array}
\right.
$$
\end{description}

\begin{description}
    \item[$(H_1)$] $\displaystyle \lim_{s\rightarrow 0} \frac{f(s)}{\left|s\right|^{N-1}} = 0$;
\end{description}

\begin{description}
  \item [$(H_2)$] There exists $\nu > N$  such that
    $$0 < \nu F(s) \leq f(s)s,\,\,\, \mbox{ for all } \, \left|s\right| > 0, $$ 
    where $F(s) = \int^{s}_{0} f(t) dt$;
\end{description}

\begin{description}
  \item [$(H_3)$] There exist $p>N$ and $C_p>0$ such that
 \[
  f(s) \geq C_p s^{p-1},\,\,\,\,\,\mbox{ for all } \, s \geq 0,
  \]
where
  \begin{equation}\label{NE}
   C_p >S_p^p \left\{\frac{p(\nu -N)}{\nu(p-N)}\left( \frac{\alpha_N}{\alpha_0}\right)^{N-1}\right\}^{ \frac{N-p}{N}},
   \end{equation}
and
   \[
   S_p := \inf_{v\in W_0^{1,N}(B_r(0)) \setminus \{0\}}\frac{||v||_{W_0^{1,N}(B_r(0))}}{|v|_{L^{p}(B_r(0))}};
   \]
\end{description}

\begin{description}
    \item[${(H_4)}$] $ \frac{f(s)}{s^{N-1}} \mbox{   is increasing in    } (0,+\infty);$
\end{description}

\begin{description}
  \item [${ (H_5)}$] There exist $\sigma \geq N $ and a constant $C>0$ such that
 \[
  f'(s)s - (N-1)f(s) \geq C s^{\sigma},\,\,\,\,\mbox{ for all }\, s \geq 0;
  \]
\end{description}

\begin{description}
		\item[${(H_6)}$] There exists $C_* > 0$ such that
    $$|f'(s)| \leq C_* e^{\alpha_0 |s|^{\frac{N}{N-1}}},\,\,\mbox{ for all } \, s \in \mathbb{R}.$$
\end{description}

In the sequel, without lost of generality, we suppose that $f(s) = 0$ in $(-\infty, 0)$, because we are looking for positive solutions.
\medskip 

One can find in the literature several studies concerning results of multiplicity linked the topology of the domain.
In \cite{BC1}, Benci \& Cerami proved that, for $2 < p<2^*$ and $\lambda$ sufficiently large, the number of solutions of problem
\begin{equation}\label{BC}
\left \{
\begin{array}{l}
-\Delta u + \lambda u = u^{p-1}, \quad \mbox{in} \quad \Omega, \\
u>0,\,\,\,\quad \mbox{in} \quad \Omega,\\
u=0,\quad  \mbox{on} \quad \partial\Omega,
\end{array}
\right.
\end{equation}
is affected by the topology of $\Omega$, where $\Omega \subset \mathbb{R}^{N}$ with $N \geq 3$, is a smooth bounded domain.  More precisely, they proved that (\ref{BC}) has at least $cat(\Omega)$   distinct solutions.
Later, in \cite{BC3}, they studied the problem
\begin{equation}\label{hand}
\left \{
\begin{array}{l}
-\epsilon^{2} \Delta u + u = f(u), \quad \mbox{in} \quad \Omega, \\
u>0,\quad \mbox{in} \quad \Omega, \\
u=0,\quad \mbox{on} \quad  \partial\Omega,
\end{array}
\right.
\end{equation}
where $\epsilon > 0$, $ \Omega\subset \RR^N \,(N\geq 3)$ and $f \in C^{1,1}(\RR^+,\RR)$ has subcritical growth. There, Benci \& Cerami used Lusternik-Schnirelman category theory to show that if $\epsilon$ is a small parameter, the problem $(\ref{hand})$ has at least $cat(\Omega)+1 $ solutions.  The reader can find more results involving Lusternik-Schnirelman category in Cerami \& Passasseo \cite{BP0}, Alves \& Ding \cite{CA4}, Rey \cite{Rey} and Bahri \& Coron \cite{BahC} and their references.

Motivated by results proved in \cite{BC1}, Alves in \cite{CA1} showed the existence of at least $cat(\Omega)$ positive solutions for the quasilinear problem   
\begin{equation} \label{CA1}
\left \{
\begin{array}{lll}
-\Delta_{p} u + |u|^{p-2}u = f(u),\quad \mbox{in} \quad  \Omega_\lambda,\\
\quad u > 0,\quad  \quad \mbox{in} \quad  \Omega_\lambda, \\
\quad u=0,\quad \quad \mbox{on} \quad  \partial\Omega_\lambda,
\end{array}
\right.
\end{equation}
where $\Omega_\la=\la\Omega$, $\la$ is a positive parameter, $2 \leq p < N$ and $f$ is a function with subcritical growth.
Succeeding this study, Alves \& Soares \cite{CS} considered the problem
\begin{equation}\label{CA2}
\begin{array}{lll}
-\epsilon^p\Delta_{p} u + (1+\la A(x))|u|^{p-2}u = f(u),\quad \mbox{in} \quad  \RR^N, \\
\end{array}
\end{equation}
where $\la,\ep >0$ are parameters and $A$ satisfying $(A_1)-(A_2)$. They proved the existence of $cat(int A^{-1}(0))$ positive solutions, for all sufficiently large $\la$ and small $\ep$.  The problem $(\ref{CA2})$ was motivated by a paper due to Bartsch \& Wang \cite{BW}, which have established the existence of at least $cat(int A^{-1}(0))$ positive solutions for the problem 
\begin{equation}\label{TB}
\begin{array}{lll}
-\Delta u + (1+\la A(x))u = u^{p-1},\quad \mbox{in} \quad \RR^N, \\
\end{array}
\end{equation}
for $N\geq 3$ and $p$ close to $2^*=\frac{2N}{N-2}$.

\medskip

The motivation of the present paper comes from \cite{CA1}, \cite{AC6}, \cite{CS} and \cite{BW}, as well as by the fact that we did not find in the literature any paper dealing with the existence of the positive solutions for the problem $(P_{\mu,\epsilon})$ involving a nonlinearity with exponential critical growth. Quasilinear  problems of the type
$$
\left \{
\begin{array}{lll}
-\Delta_{N} u = f(u),\quad \mbox{in} \quad \Omega\subset \mathbb{R}^{N}\, (N\geq 2),\\
\quad u \in W_{0}^{1,N}(\Omega),
\end{array}
\right.
$$
where $f(u)$ behaves like $exp(\alpha|u|^{\frac{N}{N-1}})$, as $\left|u\right|\rightarrow \infty$,  have been extensively analyzed by several authors, see \cite{AY, AC6, ACL, AL, JMBO, JMBO1, OMS, dooms, RP, sergio, tonkes, WYZ0} and their references. These articles were motivated by the Trudinger-Moser inequality
\[
\sup_{||u||_{W_0^{1,N}(\Omega)} \leq 1} \int_{\Omega} e^{\alpha|u|^{\frac{N}{N-1}}} dx \leq C(N,|\Omega|) , \,\,\,\,\,\mbox{ for all } \, \alpha  \leq \alpha_N = N\omega_{N-1}^{\frac{1}{N-1}}> 0,
\]
where $\omega_{N-1}$ is the  $(N-1)-$dimensional measure of the $(N-1)-$sphere.

\medskip

The our main result is the following  

\begin{thm}\label{teo 9}
Suppose that $({A_1})- ({A_2})$ and $({H_0})- ({H_6})$ hold. Then,  there exists $\epsilon^*>0$ such that for any $0 < \epsilon < \epsilon^*$,   there exists $\mu^*(\epsilon)>0$ such that the problem  $(P_{\mu,\epsilon})$ has at least $cat(\Omega)$ solutions for $\mu > \mu^*(\epsilon)$.
\end{thm}

In the proof of Theorem \ref{teo 9},  it is crucial to understand the behavior the some minimax levels of the energy functional associated with the limit problem of $(P_{\mu,\epsilon})$, given by 
$$
\left \{
\begin{array}{lll}
-\epsilon^{N}\Delta_{N} u + |u|^{N-2}u = f(u),\quad \mbox{in} \quad \Omega,\\
\quad u > 0,\quad \mbox{in} \quad   \Omega, \\
\quad u=0,\quad \mbox{on} \quad \partial\Omega,
\end{array}
\right.\eqno{ (LP)_{\ep}}
$$
where $\epsilon$ is a positive parameter, $\Omega \subset \mathbb{R}^{N} , N\geq2, $
 is a bounded smooth domain with $0 \in \Omega$. Since we did not find any result involving this study for this problem with $f$ having exponential critical growth, we were naturally taken to prove the following result:

\begin{thm}\label{teo 8}
Suppose that $f$ is a function satisfying $({H_0})- ({H_6})$. Then,  there exists a constant $\e^*>0$ such that for $0 < \e < \e^*$, the problem 
$(LP){_\ep}$ has at least $cat(\Omega)$ positive solutions. 
\end{thm}

\medskip

We would like point out that, if $Y$ is a closed subset of a topological space $X$, the Lusternik-Schnirelman category $cat_{X}(Y)$ is the least number of closed and contractible sets in $X$ which cover $Y$. Hereafter, $cat(X)$ denotes $cat_{X}(X)$.

\medskip

The paper is organized as follows. In Section 2, we state the Trudinger-
Moser inequalities and show some technical lemmas. Section 3 is devoted to study
of the problem $(LP)_{\ep}$.  Section 4 is devoted to show some technical results related to problem $(P_{\mu,\epsilon})$, while   in Section 5 we prove the Theorem  \ref{teo 9}.

\medskip

In this work we make use of the following notations:
\begin{itemize}
	\item $L^t(\Omega), \,1 \leq t<\infty$ denotes the Lebesgue spaces with the usual norm $\left|u\right|_t=\left(\displaystyle \int_{\Omega}\left|u\right|^tdx\right)^{\frac{1}{t}}$;
		\item $W_0^{1,N}(\Omega)$ denote the Sobolev space with the usual norm \linebreak $\left\|u\right\|=\left(\displaystyle \int_{\Omega} \left(\left|\nabla u \right|^{N} + |u|^N \right)dx \right)^{\frac{1}{N}}$;
		\item $C,\, C_0,\,C_1,\,C_2,...$ denote positive generic constants.
\end{itemize}

\section{Results involving exponential critical growth}

\hspace{0.5 cm}The hypothesis $(H_0)$ is motivated by the following estimates proved by Trudinger \cite{trudinger} and Moser \cite{moser}.

\begin{lem} \noindent {\bf (Trudinger-Moser inequality for bounded domains)}\label{lema 1.1}
Let $\Omega \subset \RR^N$ ($N\geq 2$) be a bounded domain. Given any  $ u \in W_0^{1,N}(\Omega)$, we have
$$
\int_{\Omega}
e^{\alpha\left|u\right|^{\frac{N}{N-1}}}dx
< \infty, \,\,\,\, \mbox{ for every }\,\,\alpha >0.
$$
Moreover, there exists a positive constant $C=C(N,|\Omega|)$ such that
\[
\sup_{||u||\leq 1} \int_{\Omega} e^{\alpha|u|^{\frac{N}{N-1}}} dx \leq C , \,\,\,\,\,\,\, \mbox{for all } \, \alpha  \leq \alpha_N,
\]
where $\alpha_N = N\omega_{N-1}^{\frac{1}{N-1}}> 0$ and $\omega_{N-1}$ is the  $(N-1)-$dimensional measure of the $(N-1)-$sphere.
\end{lem}

The next result is a version of the Trudinger-Moser inequality for whole $\mathbb{R}^{N}$, and its proof can be found in Cao \cite{Cao}, for $N = 2$, and Bezerra do \'O  \cite{JMBO}, for the case $N \geq 2$. 

\begin{lem} \noindent {\bf (Trudinger-Moser inequality for unbounded domains)}\label{lema 1.2}
Given any  $ u \in W^{1,N}(\mathbb{R}^{N})$ with $ N \geq 2$, we have
$$
\int_{\mathbb{R}^N}
\left(e^{\alpha\left|u\right|^{\frac{N}{N-1}}}-S_{N-2}(\alpha,u)\right)dx
< \infty,\,\,\,\, \mbox{ for every }\,\,\alpha >0.
$$
 Moreover, if $\left| \nabla
u\right|^{N}_{N}\leq 1,\,\left|u\right|_{N} \leq M < \infty $ and
$\alpha < \alpha_{N}$, then there
exists a positive constant $C=C(N,M,\alpha)$ such that
$$
\int_{\mathbb{R}^N} \left(e^{\alpha\left|u\right|^{\frac{N}{N-1}}}-S_{N-2}(\alpha,u)\right)dx \leq C,
$$
where
$$
S_{N-2}(\alpha,u) = \sum^{N-2}_{k=0}\frac{\alpha^{k}}{k!}\left|u\right|^{\frac{Nk}{N-1}}.
$$
\end{lem}

The Trudinger-Moser inequalities will be strongly utilized throughout
this work in order to deduce important estimates. In the sequel, we state
some technical lemmas found in \cite{AL}, which will be essential to carry out the proof of our
results.

\begin{lem}\label{lema 1.3} Let  $ \alpha > 0$ and $t >1$. Then, for every each $\beta > t$, there exists a constant $C=C(\beta,t)
> 0$ such that
$$
 \left(e^{\alpha\left|s\right|^{\frac{N}{N-1}}}-S_{N-2}(\alpha,s)\right)^{t} \leq C
 \left(e^{\beta\alpha\left|s\right|^{\frac{N}{N-1}}}-S_{N-2}(\beta\alpha,s)\right).
$$
 \end{lem}

\begin{lem} \label{alphat11} Let $(u_{n})$ be a sequence in $W^{1,N}(\mathbb{R}^{N})$ with 
$$
\limsup_{n \to +\infty} \|u_n \|^{N}  < \left(
\frac{\alpha_N}{\alpha_0}\right)^{N-1}.
$$ 
Then, there exist  $\alpha > \alpha_0$, $t > 1$ and $C > 0$  independent of $n$, such that
\[
\int_{\mathbb{R}^{N}}\left(e^{\alpha|u_n|^{\frac{N}{N-1}}} -
S_{N-2}(\alpha,u_n)\right)^t dx \leq C, \,\,\,\,\mbox{for all } \, n \geq
n_0,
\]
for some $n_0$ sufficiently large.
\end{lem}

\begin{cor} \label{lema 2.5} Let $\cal{B}$ a bounded domain in $\RR^N$ and $(u_{n})$ be a sequence in $W_{0}^{1,N}(\cal{B})$ with 
$$
\limsup_{n \to +\infty} \|u_n \|^{N} < \left(
\frac{\alpha_N}{\alpha_0}\right)^{N-1}.
$$ 
Then, there exist  $\alpha > \alpha_0$, $t > 1$ and $C > 0$  independent of $n$, such that
\[
\int_{\cal{B}}e^{t \alpha |u_n|^{\frac{N}{N-1}}} dx \leq C, \,\,\,\,\mbox{for all } \, n \geq
n_0,
\]
for some $n_0$ sufficiently large.
\end{cor}

\section{The limit problem $(LP){_\ep}$}
\hspace{0.5 cm}Using  standard arguments, we know that $(LP){_\ep}$ is equivalent to the problem 
$$
\left \{
\begin{array}{lll}
-\Delta_{N} u + |u|^{N-2}u = f(u),\quad \mbox{in} \quad \Omega_\ep,\\
\quad u > 0,\quad \mbox{in} \quad   \Omega_\ep, \\
\quad u=0,\quad \mbox{on} \quad \partial\Omega_\ep,
\end{array}
\right.\eqno{ (P_{\ep})}
$$
where $\Omega_\ep = \frac{1}{\ep}\Omega$ and  $ \ep>0$. Let $ I_{\e} : W^{1,N}_0(\Omega_\e) \to \mathbb{R} $  given by
$$
I_{\e}(u) = \frac{1}{N}\int_{\Omega_\e} \left(\left|\nabla u
\right|^{N} + |u|^N \right)dx  - \int_{\Omega_\e} F(u)dx
$$
be the functional associated with $(P_\e)$ and define the Nehari manifold
\[
{\cal M}_{\e} : = \left\{u \in W^{1,N}_0(\Omega_\e)\setminus \{0\} \, ; \,I'_{\e}(u)u =0\right\}.
\]
In what  follows, we consider $B_{\frac{r}{\e} } := B_{\frac{r}{\e} }(0)$ and denote by $ I_{\e,B} : W^{1,N}_0(B_{\frac{r}{\e}}) \to \mathbb{R} $ the functional 
$$
I_{\e,B}(u) = \frac{1}{N}\int_{B_{\frac{r}{\e}}} \left(\left|\nabla u
\right|^{N} + |u|^N \right)dx  -  \int_{B_{\frac{r}{\e} }} F(u)dx,
$$
whose corresponding Nehari manifold is given by
$$
      {\cal M}_{ \e,B} := \left\{u \in W^{1,N}_0(B_{\frac{r}{\e} })\setminus \{0 \} \, ; \, I'_{\e,B}(u)u =0 \right\}.
$$

Using well known arguments, if $c_\e$ and $b_\e$ denote the mountain pass levels associated with $I_{\e}$ and $I_{\e,B}$  respectively, then they satisfy
\[
b_\e = \inf_{u \in {\cal M}_{\e,B}} I_{\e,B}(u)
\]
and
\[
c_\e = \inf_{u \in {\cal M}_{\e}} I_{\e}(u).
\]
Apart from the above problems, we also consider the problem 
$$
\left \{
\begin{array}{lll}
-\Delta_{N} u + |u|^{N-2}u = f(u),\quad \mbox{in} \quad  \mathbb{R}^N,\\
\quad u > 0,\quad  \quad \mbox{in} \quad  \mathbb{R}^N, \\
\quad u \in W^{1,N}(\mathbb{R}^N),
\end{array}
\right.\eqno(P_\infty)
$$
whose  functional corresponding to a variational approach is \linebreak $ I_\infty :
W^{1,N}(\mathbb{R}^N) \to \mathbb{R} $ given by
$$
I_\infty(u) = \frac{1}{N}\int_{\mathbb{R}^N} \left(\left|\nabla u
\right|^{N} + |u|^N \right)dx  - \int_{\mathbb{R}^N} F(u)dx.
$$
The Nehari manifold associate to $ I_\infty$ is defined by  
$$
  {\cal M} _{\infty} : = \left\{u \in W^{1,N}(\mathbb{R}^N)\setminus \{0\} \, ; \, I'_{\infty}(u)u =0 \right\}
$$
and let us denoted by  $c_\infty$ the mountain pass level of $I_\infty$, which satisfies
\[
c_\infty = \inf_{u \in {\cal M}_\infty} I_{\infty}(u).
\]
\subsection{A result of compactness}

\hspace{0.5 cm}In this section, we will establish a result of compactness for $I_\e$ restricts to Nehari manifold ${\cal M}_{\e}$. Moreover,  we will prove that a critical point of $I_\e$ on ${\cal M}_{\e}$ is also critical point of $I_\e$ in $W_{0}^{1,N}(\Omega_\e)$. Initially, we need to study the behavior of levels $b_\e$, $c_\e$  and $c_\infty$.

\begin{lem}\label{lema 3.1} The level $b_\e $ satifies
\[
b_\e <\left(\frac{1}{N}- \frac{1}{\nu}\right)\left(\frac{\alpha_N}{\alpha_0}\right)^{N-1}, \quad \mbox{ for all } \ep \in (0,1).
\]
\end{lem}

\proof  \, 
Let  $\varphi \in W_0^{1,N}( B_r(0))\setminus \{0\}$  be a function satisfying
\[
\frac{||\varphi||}{|\varphi|_p} = S_p :=
\inf_{v\in W_0^{1,N}(B_r(0)) \setminus \{0\}}\frac{||v||}{|v|_{p}}.
\]
As $\ep \in (0,1)$, it follows that $B_r(0) \subset B_\frac{r}{\ep}(0)$. Then,  we can consider  that 
$\varphi \in W^{1,N}_0( B_\frac{r}{\ep}(0))$, and so, by definition of $b_\ep$,
\[
b_\ep \leq \max_{t \geq 0} I_\ep(t{\varphi})\leq \max_{t \geq
0}\left\{\frac{t^N}{N}||{\varphi}||^N - \int_{B_r(0)}F(t{\varphi})dx\right\}.
\]
By $(H_3)$,
\[
\frac{b_\ep}{|{\varphi}|_{p}^N} \leq \max_{t \geq
0}\left\{\frac{t^N}{N}S_p^N -  \frac{C_p}{p}t^p
|{\varphi}|_{p}^{p-N}\right\},
\]
where the maximum is attained at 
\[
 t_0 = C_p^{\frac{1}{N-p}}|{\varphi}|_{p}^{-1}S_p^{\frac{N}{p-N}}.
 \]
Consequently,
\[
b_\ep \leq \left(\frac{1}{N}- \frac{1}{p}\right)C_p^{\frac{N}{N-p}}
S_p^{\frac{pN}{p-N}},
\]
and by (\ref{NE}),  
\[
b_\ep <\left(\frac{1}{N}- \frac{1}{\nu}\right)\left(\frac{\alpha_N}{\alpha_0}\right)^{N-1}.
\]
\eproof

\begin{rmk}\label{ab0}
In order to simplify, we denote
\begin{equation}
\Lambda:=\left(\frac{1}{N}- \frac{1}{\nu}\right)\left(\frac{\alpha_N}{\alpha_0}\right)^{N-1}.
\end{equation}
Observe that, by definition
\begin{equation}
0<c_\infty \leq c_\e \leq b_\e <\Lambda, \,\,\,\,\mbox{ for all } \, \ep>0.
\end{equation}
\end{rmk}

The next proposition is an important result of compactness involving $I_\infty$.
\begin{prop}\label{prop 3.1}
Let $(u_n) \subset {\cal M}_\infty$ be a sequence satisfying
\begin{equation}\label{red}
I_\infty (u_n) \rightarrow  c_\infty.
\end{equation}
Then, $(u_n)$ admits a subsequence converging strongly in $W^{1,N}(\RR^N)$, or there exists $(y_n) \subset \mathbb{R}^N$ with $|y_n|\rightarrow \infty$ such that
\[
v_n(x) = u_n (x + y_n) \rightarrow v \,\,\, \mbox{ in } \,\,\, W^{1,N}(\RR^N),
\]
where $v \in {\cal M}_\infty$ and $I_\infty (v) = c_\infty$.
\end{prop}

\proof \, Let $(u_n) \subset {\cal M}_\infty$  verifying (\ref{red}).
 We claim that 
\begin{equation}\label{pvee}
I'_\infty(u_n)\rightarrow 0  \,\,\,\,  \mbox{ in  } \,\,\,\, (W^{1,N}(\RR^N))'.
\end{equation}
In fact, using Ekeland Variational Principle (see \cite{Ekeland}), there exists a sequence $(w_n)\subset {\cal M}_{\infty}$ verifying
\[
w_n = u_n + \o, \,\,\,\,\, I_\infty (w_n) \rightarrow c_\infty
\]
and
\begin{equation}\label{parr}
I'_\infty (w_n) - \ell_n E'_\infty (w_n) = \o,
\end{equation}
where $(\ell_n) \subset \mathbb{R}$ and $E_\infty (w) = I'_\infty (w)w$, for all $w \in W^{1,N}(\RR^N)$. By a straightforward computation, 
\[
E'_\infty (w_n)w_n = \int_{\RR^N} \left[ f'(w_n)w_n - (N-1)f(w_n)\right] w_n dx.
\]
Thereby, by $(H_5)$, there exist $ \sigma \geq N $ and a constant $ C> 0 $ such that 
\begin{equation} \label{ab33}
-E'_\infty (w_n)w_n  \geq C \int_{\RR^N} |w_n|^{\sigma +1}\,dx.
\end{equation}
Using the last expression, we can prove that there exists $\delta > 0$ such that $|E'_\infty(w_n)w_n| \geq \delta$ for all $n \in \mathbb{N}$. Indeed, suppose by contradiction that there exists a subsequence, still denoted by $(w_n)$, such that
\begin{equation} \label{ab4}
E'_\infty (w_n)w_n = \o.
\end{equation}
From (\ref{ab33}), 
\[
\int_{\RR^N} |w_n|^{\sigma +1}\,dx =\o. 
\]
Hence, by interpolation
\begin{equation}\label{verr}
\int_{\RR^N} |w_n|^{\tau}\,dx = \o, \,\,\, \mbox{for all} \,\,\,\,  \tau \geq \sigma +1.
\end{equation}
From definition of $ (w_n) $, together with $(H_2)$,  it follows that 
\begin{equation}\label{wil}
\left(\frac{1}{N} - \frac{1}{\nu}\right)||w_n||^N \leq I_\infty(w_n) - \frac{1}{\nu} I'_\infty (w_n)w_n = c_\infty +\o,
\end{equation}
showing that $(w_n)$ is bounded in $W^{1,N}(\mathbb{R}^N)$ with 
\[
\limsup_{n\rightarrow\infty} ||w_n||^N < \frac{c_\infty}{\left(\frac{1}{N}-\frac{1}{\nu}\right)} .
\]
Therefore, by Remark \ref{ab0}, 
\[
\limsup_{n\rightarrow\infty} ||w_n||^N < \left(\frac{\alpha_N}{\alpha_0}\right)^{N-1}.
\]
As $(w_n) \subset W^{1,N}(\RR^N)$, by Lemma \ref{alphat11}, there exist $\alpha > \alpha_0$, $t
> 1$ and $C > 0$, independent of $ n $, such that
\begin{equation}\label{ab22}
\int_{\mathbb{R}^{N}}\left(e^{\alpha|w_n|^{\frac{N}{N-1}}} -
S_{N-2}(\alpha,w_n)\right)^t dx \leq C, \,\,\,\,\mbox{for all } \, n \geq
n_0,
\end{equation}
for some $n_0$ sufficiently large.

From $({H}_0)$ and $({H}_1)$, for each  $ \eta> 0 $ and $ s \geq 1 $, there exists $ C> $ 0 such that
\begin{equation}
||w_n||^N  \leq  \eta |w_n|^N_N + C\int_{\RR^N}|w_n|^s \left(e^{\alpha|w_n|^{\frac{N}{N-1}}} -
S_{N-2}(\alpha,w_n)\right)dx.
\end{equation}
Choosing $ \eta $ small enough, $s \geq \sigma +1$ and using  H\"{o}lder's inequality  together with (\ref{ab22}), we find
\begin{equation}\label{funny}
||w_n||^N  \leq  C|w_n|_{st_1}^s,
 \end{equation}
where $t_1=\frac{t}{t-1}$.
Therefore, from (\ref{verr}), 
\[
||w_n||^N = \o,
\]
that is, $w_n \rightarrow 0$ in $W^{1,N}(\RR^N)$. On the other hand, from (\ref{funny}), 
\[
||w_n||^{s-N}\geq C_2 > 0,
\]
for some positive constant $C_2$, which is a contradiction. This contradiction yields there exists $ \delta> 0 $ such that
\begin{equation}\label{onee}
|E'_\infty(w_n)w_n| \geq \delta, \,\,\,\,\, \mbox{ for all } \, n \in \NN.
\end{equation}
Now, by (\ref{parr}) 
\[
\ell_n E'_\infty(w_n)w_n = \o,
\]
and so, $\ell_n = \o$. Since $(w_n)$ is bounded, it is not difficult to prove that $ (E'_\infty(w_n))$ is bounded. Using again (\ref{parr}), we can ensure that 
\begin{equation}\label{pveee}
I'_\infty(w_n)\rightarrow 0 \,\,\, \mbox{  in  } \,\,\, \left(W^{1,N}(\RR^N)\right)'.
\end{equation}
Therefore, without loss generality, we can assume that
\begin{equation}\label{pvee}
I'_\infty(u_n)\rightarrow 0  \,\,\,\,  \mbox{ in  } \,\,\,\, (W^{1,N}(\RR^N))'.
\end{equation}
As in (\ref{wil}), $(u_n)$ is bounded in $W^{1,N}(\mathbb{R}^N)$ and
\[
\limsup_{n\rightarrow \infty}||u_n||^N   < \left(\frac{\alpha_N}{\alpha_0}\right)^{N-1}.
\]
Thus, there exists $u \in W^{1,N}(\mathbb{R}^N)$ such that, for a subsequence, $u_n \rightharpoonup  u$ in $W^{1,N}(\mathbb{R}^N)$. Now, we divide the proof into two cases:

\noindent {\bf Case I: }  $ u \neq 0 $.
Repeating the same arguments employed in \cite[Lemma 3]{AC6}, we see that 
\begin{description}
    \item [$(i)$] $\nabla u_n(x) \rightarrow \nabla u(x) \,\,\, \mbox{  a.e.   in  } \,\,\, \mathbb{R}^N$;
    \item [$(ii)$] $ \int_{\mathbb{R}^N}|\nabla u_n|^{N-2}\nabla u_n \nabla v dx \rightarrow  \int_{\mathbb{R}^N}|\nabla u|^{N-2}\nabla u \nabla v dx, \,\,\, \forall \, v \in W^{1,N}(\mathbb{R}^N)$;
    \item [$(iii)$] $\int_{\mathbb{R}^N}| u_n|^{N-2} u_n v dx \rightarrow  \int_{\mathbb{R}^N}| u|^{N-2} u  v dx, \,\,\, \forall v \,\in W^{1,N}(\mathbb{R}^N)$;
		\item [$(iv)$] $\int_{\mathbb{R}^N} f(u_n) \varphi dx \rightarrow  \int_{\mathbb{R}^N} f(u) \varphi dx, \,\,\, \forall \,  \varphi \in C_0^\infty(\mathbb{R}^N)$.
\end{description}
Using the above limits we can easily show that $I'_\infty(u)=0$. This fact together with Fatou's Lemma leads to
 \begin{eqnarray*}
 c_\infty & \leq & I_\infty (u) = I_\infty(u) -\frac{1}{\nu} I'_\infty(u)u\\
 & = & \left(\frac{1}{N} - \frac{1}{\nu}\right)||u||^N + \int_{\mathbb{R}^N} \left[\frac{1}{\nu} f(u)u - F(u)\right]\, dx \\
 & \leq & \liminf_{n \rightarrow \infty}  \left\{\left(\frac{1}{N} - \frac{1}{\nu}\right)||u_n||^N + \int_{\mathbb{R}^N}\left[\frac{1}{\nu} f(u_n)u_n - F(u_n)\right] \,dx \right\}\\
 & = & \liminf_{n \rightarrow \infty} \{ I_\infty(u_n) -\frac{1}{\nu} I'_\infty(u_n)u_n \} = c_\infty,
 \end{eqnarray*}
from where it follows that  
\[
\lim_{n\rightarrow \infty} ||u_n||^N = \lim_{n \rightarrow \infty} ||u||^N.
\]
The last limit and $(i)$ combine to give $u_n \rightarrow u$ in $W^{1,N}(\mathbb{R}^N)$.

 \medskip

\noindent {\bf Case II: }  $ u = 0 $.   First, let us claim that there exist
 $R, \kappa > 0$ and $(y_n) \subset \mathbb{R}^N$ such that
\begin{equation}\label{ab6}
\limsup_{n \rightarrow \infty} \int_{B_R(y_n)}|u_n|^N \,dx\geq \kappa.
\end{equation}
In fact, otherwise
\[
\limsup_{n \rightarrow \infty} \left[\sup_{y \in \mathbb{R}^N} \int_{B_R(y)}|u_n|^N \,dx\right]= 0
\]
and by Lions (see \cite{Lio}),
\[
u_n \rightarrow 0 \,\, \mbox{  in  } \,\, L^q(\mathbb{R}^N), \,\, \mbox{   for all  }\,\,\, q \in (N,+\infty).
\]
From $(H_0)$ and $(H_1)$, given $s > N$ and $\alpha > \alpha_0$, there exists a constant
$C = C(s) > 0$ such that 
\[
||u_n||^N = \int_{\mathbb{R}^N} f(u_n)u_n\,dx \leq \frac{1}{2}||u_n||^N + C\int_{\RR^N}|u_n|^s \left[e^{\alpha|u_n|^{\frac{N}{N-1}}}-S_{N-2}(\alpha, u_n)\right]\,dx.
\]
Combining H\"{o}lder's inequality with Trudinger-Moser inequality (Lemma \ref{alphat11}), we derive the inequality 
\[
||u_n||^N  \leq  C|u_n|_{ts}^s.
\]
Then  $||u_n||\rightarrow 0$, and so, 
\[
I_\infty (u_n) \rightarrow 0,
\]
which is a contradiction, because $I_\infty (u_n) \rightarrow c_\infty>0$. This way, (\ref{ab6}) holds.

\medskip
We can notice that $|y_n|\rightarrow +\infty$. Indeed, if $(y_n)$  is bounded, that is, $|y_n| \leq K$ for all $n \in \mathbb{N}$, we have
\[
\int_{B_{R+K}(0)} |u|^N dx  = \limsup_{n\rightarrow \infty}\int_{B_{R+K}(0)} |u_n|^N dx \geq \limsup_{n\rightarrow \infty}\int_{B_{R}(y_n)} |u_n|^N dx \geq \kappa > 0,
\]
obtaining a new contradiction, because $u=0$.

Now, we define $v_n(x) = u_n (x+y_n)$. Making change of variable, it follows that
\[
I_\infty(v_n)  = I_\infty(u_n) \,\,\,\, \mbox{ and } \,\,\,\, I'_\infty(v_n) =  I'_\infty(u_n), \,\,\,\, \forall n \in \mathbb{N},
\]
that is,
\begin{equation}\label{jam}
I_\infty(v_n)\rightarrow c_\infty \,\,\,\, \mbox{ and } \,\,\,\, I'_\infty(v_n) \rightarrow 0, \,\,\,\, \mbox{ as } n\rightarrow +\infty.
\end{equation}

\medskip

It is easy to see that, $(v_n)$ is bounded in $W^{1,N}(\RR^N)$ and there exists $v \in W^{1,N}(\RR^N) $ such that, for a subsequence, $v_n \rightharpoonup v$ in $W^{1,N}(\RR^N)$ and
\[
\int_{B_{R}(0)} |v(z)|^N dz  = \limsup_{n\rightarrow \infty}\int_{B_{R}(0)} |v_n(z)|^N dz = \limsup_{n\rightarrow \infty}\int_{B_{R}(y_n)} |u_n(x)|^N dx \geq \kappa ,
\]
showing that $v \neq 0$. Arguing as in Case I,  $v_n \rightarrow v$ in $W^{1,N}(\RR^N)$. From (\ref{jam}), it follows that $v \in {\cal M}_\infty$ and $I_\infty (v) = c_\infty$, completing the proof.
\eproof

\medskip

The next result shows that $I_\ep$ restricts to Nehari Manifold verifies the $(PS)$ condition at some levels. Hereafter, the  norm of the derivative of the restriction of $I_\e$ to ${\cal M}_\e$  at $v$ is defined as
\[
||I'_\e(v)||_*:= \sup_{w \in T_v {\cal M}_\e,\,||w||=1} I'_\e (v) w.
\]
\begin{lem} \label{lema 3.2}
The functional $I_\e$ restricts to Nehari manifold ${\cal M}_{\e}$ satisfies the $(PS)_c$ condition for all $c < \Lambda$. 
\end{lem}

\proof \, Let $(u_n) \subset {\cal M}_\e$ such that
\[
I_\e(u_n) \rightarrow c \,\,\, \mbox{  and   } \,\,\, ||I'_\e(u_n)||_*\rightarrow 0,
\]
where $c <\Lambda$. We will prove that $(u_n)$ admits a subsequence strongly convergent in  $W^{1,N}_0(\Omega_\e)$. 

Applying  the Ekeland's Variational Principle and arguing as in the proof of (\ref{pveee}), we can assume that 
$I'_\e(u_n) \rightarrow 0$. Thus, $(u_n)$  is a $(PS)_c$ sequence for  $I_\e$, and analogously to $(\ref{wil})$, we obtain that $(u_n)$ is bounded in $W^{1,N}_0(\Omega_\e)$ and
 \begin{equation}\label{per}
 \limsup_{n \rightarrow \infty} ||u_n||^N \leq \frac{c}{\left(\frac{1}{N}- \frac{1}{\nu}\right)} < \left(\frac{\alpha_N}{\alpha_0}\right)^{N-1}.
 \end{equation}
Therefore, there exists
$u \in W^{1,N}_0(\Omega_\e)$ such that, for a subsequence, $u_n \rightharpoonup u$ in $W^{1,N}_0(\Omega_\e)$.  Using standard arguments, it is easy to prove that $u$ is a critical point of $I_\e$. 

Combining (\ref{per}), the Lemma \ref{lema 2.5} and a Br\'ezis-Lieb Lemma, we have that, for $s > N$ and $\alpha>\alpha_0\,(\alpha\approx\alpha_0)$, 
\[
\int_{\Omega_\e}|u_n|^s e^{\alpha|u_n|^{\frac{N}{N-1}}}dx \rightarrow \int_{\Omega_\e}|u|^s e^{\alpha|u|^{\frac{N}{N-1}}} dx.
\]
Since
\[
|f(u_n)u_n| \leq |u_n|^N + C|u_n|^s e^{\alpha|u_n|^{\frac{N}{N-1}}},
\]
the generalized Lebesgue dominated convergence theorem ensures that
\[
\lim_{n \rightarrow \infty } \int_{\Omega_\e} f(u_n)u_n\,dx = \int_{\Omega_\e} f(u)u\,dx.
\]
Thus,
\begin{eqnarray*}
||u_n||^N & = & I'_\e(u_n)u_n + \int_{\Omega_\e} f(u_n)u_n \,dx= \int_{\Omega_\e} f(u)u \,dx+ \o\\
& = & I'_\e(u)u + ||u||^N + \o = ||u||^N + \o.
\end{eqnarray*}
Hence, $u_n \rightarrow u$ in $W^{1,N}_0(\Omega_\e)$ and $u \in {\cal M}_\e$.
\eproof

\begin{lem} \label{lema 3.3}
If $u \in {\cal M}_\e$ is a critical point of $I_\e$ in ${\cal M}_\e$, then $u $ is a nontrivial critical point of $I_\e$ in  $W^{1,N}_0(\Omega_\e)$.
\end{lem}
\proof \,
Since by hypothesis $u \in {\cal M}_\ep$ is a critical point of $I_\e$ in ${\cal M}_\e$, then $u \neq 0$ and  there exists  $\ell \in \mathbb{R}$ such that
\[
I'_\ep(u) = \ell E'_\ep(u),
\]
where $E_\ep(u) = I'_\ep(u) u $. Using the equality $I'_\ep(u)u = 0$, we also have $\ell E'_\ep(u)u =0$.  Now, proceeding as in the proof of (\ref{onee}), we know that there exists  $\delta > 0$ such that
\[
E'_\ep(u)u \leq - \delta, \,\,\,\, \mbox{for all} \,\, u \in {\cal M}_\ep,
\]
implying that $\ell = 0$,  and consequently, $I'_\ep(u) = 0$.

\eproof

\subsection{Properties of the minimax levels }

\hspace{0.6 cm} In this Section, we will prove some properties of the levels $c_\infty$, $c_\e$ and $b_\e$. For this, we need to introduce some notations. For each  $x \in \mathbb{R}^N$ and $R > r > 0$, we fix 
\[
A_{R,r,x}:=B_R(x) \setminus \overline{B}_r(x),
\]
the functional $ J_{\e,x}: W_0^{1,N}(A_{\frac{R}{\e} ,\frac{r}{\e},x})\rightarrow \RR$ by
$$
J_{\e,x}(u) = \frac{1}{N}\int_{A_{\frac{R}{\e} ,\frac{r}{\e},x}} \left(\left|\nabla u
\right|^{N} + |u|^N \right) dx -  \int_{A_{\frac{R}{\e} ,\frac{r}{\e},x}} F(u)dx
$$
and 
\[
{\cal M}_{\e,x} = \left\{u \in W^{1,N}_0(A_{\frac{R}{\e} ,\frac{r}{\e},x})\setminus \{0 \}\, ; \, J'_{\e,x}(u)u =0\right\},
\]
the  Nehari manifold related to $J_{\e,x}$. 

For each $u \in W_0^{1,N}(\Omega_\ep)\setminus \{0\}$, we set 
\[
\beta(u) := \frac{\int_{\Omega_\ep}x|\nabla u|^N dx}{\int_{\Omega_\ep}|\nabla u|^N dx}
\]
and
\[
a(R,r,\e,x) := \inf\left\{J_{\e ,x}(u)\, ; \, \beta (u) = x,\, u \in {\cal M}_{\e,x} \right\}.
\]
For $x=0$, we write simply
\[
J_\e := J_{\e,0}, \,\,\,\,A_{\frac{R}{\e} ,\frac{r}{\e}}:= A_{\frac{R}{\e} ,\frac{r}{\e},0}, \,\,\,\, {\cal \hat{M}}_{\e}:={\cal M}_{\e,0}
\]
and
\[
a(R,r,\e):=a(R,r,\e,0).
\]

The next three lemmas are crucial in our arguments, but we will omit their proof, because they follow using the same ideas found in \cite{CA1}.

\begin{lem} \label{lema 4.4} The number $a(R,r,\e)$ satisfies
\[
\liminf_{\e \rightarrow 0^+} a(R,r,\e) > c_\infty.
\]
\end{lem}

\begin{lem} \label{lema 4.5} The minimax levels verify the limits
$$
\lim_{\e \rightarrow 0^+} c_\e = \lim_{\e \rightarrow 0^+} b_\e=c_\infty.
$$
\end{lem}

\begin{lem}\label {lema 4.6} There exists $\e^* > 0$ such that for any $u \in {\cal M}_\e$ satisfying $I_\e (u) \leq b_\e$, we have 
$$
\beta (u) \in  \Omega_{\e}^+, \,\,\,\,\, \mbox{for all } \, 0 < \e < \e^*.
$$
\end{lem}

\subsection{Proof of Theorem \ref{teo 8}}

In what follows, we denote by $ u_{\e, r} \in W^{1,N}_0(B_{\frac{r}{\e}}(0))$ be a positive radial ground state solution for the functional $I_{\e,B}$, that is,
\[
I_{\e,B}(u_{\e, r}) = b_\e = \inf_{{u \in \cal M}_{\e,B} }I_{\e,B} (u) \,\,\,\, \mbox{ and } \,\,\, I'_{\e,B} (u_{\e, r}) =0.
\]
Using the function $u_{\e, r}$, we define the operator  $\Psi_r :  \Omega_{\e}^{-} \longrightarrow W^{1,N}_0(\Omega_\e)$ by
\begin{eqnarray}
\Psi_r (y)(x)=
\left \{
\begin{array}{lll}\label{1.1}
 u_{\e, r}(|x-y|),\quad x\in B_{\frac{r}{\e}}(y),\\
0,\quad  x\in \Omega^{-}_\e \setminus B_{\frac{r}{\e}}(y),
\end{array}
\right.
\end{eqnarray}\label{G1}
which is continuous and satisfies
\begin{equation}\label{G0}
\beta (\Psi_r(y)) = y, \,\,\,\,\,\forall \,y \in \Omega_{\e}^{-} .
\end{equation}
 
Using the above information, we are ready to prove the following claim

\begin{clm} \label{categ} For $0< \e <\e^*$, 
\[
cat(I_\e^{b_\e}) \geq cat(\Omega),
\]
where $I_\e^{b_\e}:=\left\{u \in {\cal M}_\e \, ;\, I_\e(u) \leq b_\e\right\}$ and $\e^*$ is given in Lemma \ref{lema 4.6}.
\end{clm}
\noindent Indeed, assume that
\[
I_\e^{b_\e} = F_1 \cup F_2 \cup \cdots\cup F_n,
\]
where  $F_j$ is closed and contractible in $I_\e^{b_\e}$, for each $j = 1,2,...,n $, that is, there exist $h_j \in C([0,1]\times F_j , I_\e^{b_\e})$ and $w_j \in F_j$ such that
\[
h_j(0,u) = u   \,\,\,\, \mbox{  and  }  \,\,\, h_j (1, u) = w_j,
\]
for all $u \in F_j$. Considering the closed sets $B_j : = \Psi_r^{-1}(F_j)$, $1 \leq j \leq n$, it follows that
\[
\Omega_{\e}^{-} = \Psi_r^{-1} (I_\e^{b_\e}) = B_1 \cup B_2 \cup \cdots \cup B_n,
\]
and defining the deformation  $g_j : [0,1] \times B_j \rightarrow \Omega_{\e}^{+}$ given by
\[
g_i(t,y) = \beta (h_j(t,\Psi_r(y))),
\]
we conclude that, by Lemma \ref{lema 4.6},  $g_i$ is well defined and thus, $B_j$ is contractible  in $\Omega_{\e}^+$ for each $j= 1,2,...,n$.
Therefore,
\[
cat(\Omega)=cat(\Omega_\e)= cat_{\Omega_{\e}^{+}}(\Omega_{\e}^{-})\leq n.
\]
finishing the proof of the claim.

Since $I_\e$ satisfies the $(PS)_c$ condition on ${\cal M}_\e$ for $c < b_\e$ (see Lemmas \ref{lema 3.2} and \ref{lema 3.1}), we can apply the Lusternik-Schnirelman category theory and the Claim \ref{categ} to ensure that  $I_\e$ has at least $cat(\Omega)$ critical points on ${\cal M}_\e$. Consequently, $I_\e$ has at least $cat(\Omega)$ critical points in $W_0^{1,N}(\Omega_\e)$ (see Lemma \ref{lema 3.3}). By maximum principle, all solutions obtained are positive.
\eproof

\section{Results of existence for the problem ($P_{\mu, \epsilon}$)}

This section is concerned with the multiplicity of solutions of the problem $(P_{\mu, \epsilon})$. First of all, we observe that the problem($P_{\mu, \epsilon}$) is equivalent to the following problem  
$$
\left \{
\begin{array}{lll}
-\Delta_{N} u + (1+\mu A(\epsilon x))|u|^{N-2}u = f(u),\quad \mbox{ in} \quad \RR^N, \\
\quad u > 0,\quad \mbox{ in} \quad \RR^N. \\
\end{array}
\right.\eqno{ (\bar{P}_{\mu,\epsilon})}
$$
For each $\epsilon, \mu > 0$, we define the Banach space $E_{\mu,\epsilon}= \left(E_{\epsilon}, \left\|\cdot\right\|_{\mu,\epsilon}\right)$, where 
$$
E_{\epsilon}=\left\{u \in  W^{1,N} \left( \mathbb{R}^N \right)\,;\, \int_{\RR^N}A(\epsilon x)|u|^{N}dx < \infty\right\}
$$
and
$$
\left\|u\right\|_{\mu,\epsilon}=\left(\int_{\RR^N} \left[|\nabla u|^N + (\mu A(\epsilon x)+1)|u|^{N}\right]dx\right)^{\frac{1}{N}}.
$$
Note that the space $E_{\mu,\epsilon}$ is continuously embedded in $W^{1,N} \left( \mathbb{R}^N \right)$. 

In what follows, we denote by $ I_{\mu, \epsilon} : E_{\mu,\epsilon} \to \mathbb{R} $  the energy functional related to $(\bar{P}_{\mu,\epsilon})$, given by
$$
I_{\mu, \epsilon}(u) = \frac{1}{N}\left\|u\right\|^{N}_{\mu,\epsilon} - \int_{\RR^N} F(u)dx.
$$
Using standard arguments,  it is possible to prove that $I_{\mu, \epsilon} \in C^{1}( E_{\mu,\epsilon},\RR)$. \\

Using the above notations, we are able to study some properties of the functional $I_{\mu, \epsilon} $.

\subsection{The Palais-Smale condition}

\hspace{0.6 cm}Throughout this section, $(u_n) \subset E_{\mu,\epsilon} $ denote a $(PS)_c$ sequence for $I_{\mu, \epsilon}$, that is, 
$$
I_{\mu, \epsilon}(u_n) \rightarrow c \,\,\,\,\mbox{ and } \,\,\,\,\, I'_{\mu, \epsilon}(u_n) \rightarrow 0.
$$

In order to prove a compactness result for $I_{\mu, \epsilon}$, we need of some lemmas.

\begin{lem} \label{lema lim}The following properties occur:

\begin{description}
	\item[i.]  $\displaystyle \limsup_{n\rightarrow \infty} \left\|u_n\right\|^{N}_{\mu,\epsilon}\leq \frac{c\nu N}{\nu -N }$;
	\item[ii.]  $c\geq 0$;
	\item[iii.] if $c=0$, then $u_n\rightarrow 0$ in $E_{\mu,\epsilon}$.
\end{description}

\end{lem}
\proof \, This lemma is an immediate consequence of the inequalities below
$$
\o +c+\o \left\|u_n \right\|_{\mu,\epsilon} \geq  I_{\mu, \epsilon}(u_n)-\frac{1}{\nu}I'_{\mu, \epsilon}(u_n)u_n \geq \left(\frac{1}{N}-\frac{1}{\nu}\right)  \left\|u_n\right\|^{N}_{\mu,\epsilon}\geq 0.
$$
\eproof

\begin{lem} \label{kings}Let $c \in (0, \Lambda)$. Then there exist $\delta >0 $ and $s> N$, independent of $\mu$ and $\epsilon$, such that
$$
\liminf_{n \rightarrow \infty} \int_{\RR^N} |u_n|^{s}dx \geq \delta.
$$
\end{lem}
\proof \, From $({H}_0)$ and $({H}_1)$, for each  $ \eta> 0 $ and $\alpha > \alpha_0$, there exists a constant $ C=C(\eta,N)> $ 0 such that
\begin{equation}\label{ple}
 |f(u_n)u_n| \leq  \eta |u_n|^N + C |u_n|^N \left(e^{\alpha|u_n|^{\frac{N}{N-1}}}-S_{N-2}(\alpha,u_n)\right).
\end{equation}
Once $c< \Lambda$, the Lemma \ref{lema lim} gives 
$$
\limsup_{n \to +\infty} \|u_n \|^{N}  < \left(
\frac{\alpha_N}{\alpha_0}\right)^{N-1}.
$$ 
Thus, by Lemma \ref{alphat11}, there exist  $\alpha > \alpha_0$, $t_1 > 1$ and $C > 0$  independent of $n$, such that for some $n_0$ sufficiently large
\begin{equation}\label{TM}
\int_{\mathbb{R}^{N}}\left(e^{\alpha|u_n|^{\frac{N}{N-1}}} -
S_{N-2}(\alpha,u_n)\right)^{t_1}dx \leq C, \,\,\,\,\mbox{for all } \, n \geq
n_0.
\end{equation}
Using H\"older's inequality together with (\ref{ple}) and (\ref{TM}), we derive that
\begin{eqnarray*}
\left\|u_n\right\|^{N}_{\mu,\epsilon} &=& I'_{\mu, \epsilon}(u_n)u_n+ \int_{\RR^{N}}f(u_n)u_ndx\\
                                        & \leq & \eta |u_n|_{N}^{N} + C \int_{\RR^{N}}|u_n|^N \left(e^{\alpha|u_n|^{\frac{N}{N-1}}}-S_{N-2}(\alpha,u_n)\right)dx  +\o \\
																				&\leq & \eta |u_n|_{N}^{N} + C |u_n|^{N}_{Nt_2} +\o,
\end{eqnarray*}
where $t_2=\frac{t_1}{t_1-1}$. Choosing $\eta$ small enough, we find
\begin{equation}\label{lu1}
\left\|u_n\right\|^{N}_{\mu,\epsilon} \leq C |u_n|^{N}_{Nt_2} +\o.
\end{equation}
On the other hand
\begin{equation}\label{lu2}
\left\|u_n\right\|^{N}_{\mu,\epsilon} = NI_{\mu, \epsilon}(u_n)+ N\int_{\RR^{N}}F(u_n)dx \geq  NI_{\mu, \epsilon}(u_n) =Nc+\o.
\end{equation}
From (\ref{lu1}) and (\ref{lu2}),
$$
0< (Nc+\o )^{t_2}\leq (2C |u_n|^{N}_{Nt_2})^{t_2}.
$$
Choosing $s=Nt_2$, it follows that
$$
\liminf_{n \rightarrow \infty} \int_{\RR^N} |u_n|^{s}dx \geq \delta,
$$
with $\delta$ and $s$ independent of $\mu$ and $\epsilon$, finishing the proof. \eproof

\begin{lem}\label{jovi} Given $\epsilon >0$, $s>N$ and $\eta >0$, there exist $M_{\eta} > 0$, independent of $\epsilon$,  and $R_{\eta} >0$, such that 
$$
\limsup_{n \rightarrow \infty} \int_{B^{c}_{R_{\eta}}} |u_n|^{s}dx < \eta, \quad \forall \mu \geq M_{\eta}.
$$

\end{lem}
\proof \,
For $R>0$, fix
$$
X_R=\left\{x \in \RR^N \, ; \, |x| > R , \, A(\epsilon x) \geq M_0\right\}
$$
and
$$
Y_R=\left\{x \in \RR^N \, ; \, |x| > R , \, A(\epsilon x) < M_0\right\},
$$
where $M_0$ is given $(A_2)$.
Observe that, 
\begin{equation}\label{lu3}
\int_{X_R}|u_n|^Ndx \leq \frac{1}{\mu M_0 +1}\int_{X_R}(\mu A(\epsilon x) +1)|u_n|^N dx\leq\frac{||u_n||^{N}_{\mu, \epsilon}}{\mu M_0 +1}
\end{equation}
and \begin{equation}\label{lu4}
\int_{Y_R}|u_n|^Ndx \leq \left(\int_{Y_R}|u_n|^s dx\right)^{\frac{N}{s}} [m (Y_R)]^{\frac{1}{s_1}} \leq C ||u_n||^{N}_{\mu, \epsilon} [m (Y_R)]^{\frac{1}{s_1}}.
\end{equation}
 Using interpolation inequality for $N<s<q$, we can infer that
\begin{equation}\label{lu5}
|u_n|_{L^s(B^{c}_R)} \leq |u_n|^{\theta}_{L^N( B^{c}_R)} |u_n|^{1-\theta}_{L^q( B^{c}_R)} \leq |u_n|^{\theta}_{L^N( B^{c}_R)} ||u_n||^{1-\theta}_{\mu, \epsilon},
\end{equation}
for some $\theta \in (0,1)$.
From (\ref{lu3})-(\ref{lu5}) and Lemma \ref{lema lim}, there exists $K>0$ such that
\begin{equation}\label{lu6}
\limsup_{n\rightarrow\infty}|u_n|_{L^s( B^{c}_R)}  \leq K \left( \frac{1}{\mu M_0 +1} + [m (Y_R)]^{\frac{1}{s_1}}\right)^{\frac{\theta}{N}}  .   
\end{equation}
From $(A_2)$,
$$
\lim_{R\rightarrow \infty}m (Y_R))=0
$$
which together with (\ref{lu6}) implies that, given $\eta >0$, we can fix $R_{\eta}>0$ and $M_{\eta}>0$ such that
\begin{equation*}
\limsup_{n\rightarrow\infty}|u_n|_{L^s(B^{c}_{R_{\eta}})}  < \eta, \quad \mu \geq M_{\eta}.
\end{equation*}
\eproof

\begin{prop} \label{fund} Given $\epsilon >0$,  there exists $\bar{\mu}$, independent of $\epsilon$, such that $I_{\mu, \epsilon}$ satisfies the $(PS)_c$ condition, for all  $\mu \geq \bar{\mu}$ and $ 0 < c < \Lambda$. Moreover, the limit of any $(PS)_c$ sequence is nontrivial.
\end{prop} 

\proof \, Let $(u_n) \subset E_{\mu, \epsilon}$ be a $(PS)_c$ sequence for $I_{\mu, \epsilon}$ with $0 < c < \Lambda$.
Since $(u_n)$ is bounded (see Lemma \ref{lema lim}), there exists $u \in  E_{\mu, \epsilon}$ such that, for some subsequence, 
\[
\left\{
\begin{array}{cl}
u_n \rightharpoonup u& \mbox{ in }   E_{\mu, \epsilon},\\
u_n(x) \rightarrow  u(x) &\mbox{ a.e. in } \RR^N,\\
 u_n \rightarrow  u&  \mbox{ in }  L_{loc}^{t}(\RR^N)
\mbox{  for  }  t \geq 1.
\end{array}\right.
\]
Setting $v_n:=u_n-u$, using the above limits and following the methods used in \cite{AC6}, we know that 
\begin{description}
    \item [$(a)$]$I_{\mu,\epsilon}(v_n) = I_{\mu,\epsilon}(u_n) - I_{\mu,\epsilon}(u) + \o$;
    \item [$(b)$]$I'_{\mu,\epsilon}(v_n) =  \o$;
    \item [$(c)$] $I'_{\mu,\epsilon}(u) =  0$.
 \end{description}
 
Let $c'=c - I_{\mu,\epsilon}(u)$. From $(a)$ and $(b)$, $(v_n)$ is a $(PS)_{c'}$ sequence for $I_{\mu, \epsilon}$ and
$$
c'=c - I_{\mu,\epsilon}(u)=c-I_{\mu,\epsilon}(u)+\frac{1}{\nu}I'_{\mu,\epsilon}(u)u \leq c-\left(\frac{1}{N}-\frac{1}{\nu}\right)||u||^{N}_{\mu, \epsilon}\leq c.
$$
Furthermore,  we claim that $c' =0$. Indeed, by Lemma \ref{lema lim}, $c'\geq 0$. Supposing by contradiction that $c' > 0$, the Lemma \ref{kings} guarantees that there exist $\delta >0 $, $s> N$, independent of $\mu$ and $\epsilon$, such that
\begin{equation} \label{rit}
\liminf_{n \rightarrow \infty} \int_{\RR^N} |v_n|^{s}dx \geq \delta. 
\end{equation}
By Lemma \ref{jovi}, with $\eta=\frac{\delta}{2}$, there exist $\bar{\mu}$ (independent of $\epsilon$)  and $R$, such that 
\begin{equation}\label{mil}
\limsup_{n \rightarrow \infty} \int_{B^{c}_{R}} |v_n|^{s}dx < \frac{\delta}{2},
\end{equation}
for all $\mu \geq \bar{\mu}$. Hence,
\begin{equation} \label{rit1}
\delta \leq   \liminf_{n \rightarrow \infty} \int_{ B^{c}_{R}} |v_n|^{s}dx + \liminf_{n \rightarrow \infty} \int_{ B_{R}} |v_n|^{s}dx <  \frac{\delta}{2},
\end{equation}
which is an absurd, and so, we must have $c' =0$. Thereby,  by Lemma \ref{lema lim}, we have that  $v_n \rightarrow 0$, that is,  $u_n \rightarrow u$. As $c>0$, it follows that $u \not=0$. \eproof

\medskip

\begin{prop}  \label{rox}Let $0< \ep <1 $  fixed and $(u_n) \subset W^{1,N}(\RR^N)$  be a  sequence of solutions of $(P_{\mu_n, \ep})$, with $\mu_{n} \rightarrow \infty$ and $ \displaystyle \limsup_{n \rightarrow \infty} I_{\mu_n, \epsilon}(u_n) < \Lambda$. Then,  there exists $u_{\ep} \in   W^{1,N}(\RR^N)$ solution of $(P_{\ep})$ such that, for some subsequence,  

\begin{description}
	\item[i)] $u_{\ep} \in W^{1,N}_0(\Omega_{\ep})$;
	\item[ii)] $u_n \rightarrow u_{\ep}$ strongly in  $ W^{1,N}(\RR^N)$;
	\item[iii)] $\mu_n \int_{\RR^N}A(\ep x)|u_n|^{N}dx\rightarrow 0$;
	\item[iv)] $||u_n-u_{\ep}||^{N}_{\mu_n, \epsilon}\rightarrow 0$. 
\end{description}
\end{prop} 

\proof \, First of all, we observe that $(||u_n||_{\mu_n, \epsilon})$ is bounded in $\RR$ with 
$$
\limsup_{n \to \infty}\|u_n\|_{\mu_n,\epsilon} < \left(
\frac{\alpha_N}{\alpha_0}\right)^{N-1}.
$$
Thereby, there exists $u_{\ep} \in  W^{1,N}(\RR^N)$ such that, for a subsequence, 
\[
\left\{
\begin{array}{cl}
u_n \rightharpoonup u_{\ep}& \mbox{ in }   W^{1,N}(\RR^N),\\
u_n(x) \rightarrow  u_{\ep}(x) &\mbox{ a.e. in } \RR^N,\\
 u_n \rightarrow  u_{\ep}&  \mbox{ in }  L_{loc}^{t}(\RR^N)
\mbox{  for  }  t \geq 1.
\end{array}\right.
\]
For $k \in \NN$, we define
\[
C_k=\left\{ x \in \RR^N \,; \, A_\e (x) \geq \frac{1}{k}\right\}, \,\,\, \mbox { where }\,\, A_\e (x) := A (\epsilon x).
\]
Note that
\begin{equation*}
\int_{C_k}|u_n|^{N}dx\leq  \int_{C_k}k A(\epsilon x) |u_n|^{N} dx\leq \frac{k}{\mu_n}||u_n||^{N}_{\mu_n, \epsilon}= \o.
\end{equation*}
By Fatou's Lemma
\begin{equation*}
\int_{C_k}|u_\ep|^{N}dx=0,
\end{equation*}
and consequently, $u_{\ep}=0$ a.e. in $C_k$. Since $\RR^N \setminus A_{\e}^{-1}(0)= \cup^{\infty}_{k=1} C_k$, and $A_{\e}^{-1}(0)=\bar{\Omega}_{\e}\cup D_\e$ with  $m(D_\e)=0$, it follows that $u_\e=0$ a.e. in $\RR^N \setminus \bar{\Omega}_{\e} $. As $\partial\Omega_\e$ is smooth, we conclude that $u_{\ep} \in W^{1,N}_0(\Omega_{\ep})$.

Repeating the same arguments employed in the proof of Lemma \ref{jovi}, we know that
given any  $\eta >0$ and $s>N$, there exist $M_{\eta}$ and $R_{\eta}$,  such that 
\begin{equation}\label{fir}
\limsup_{n \rightarrow \infty} \int_{B^{c}_{R_{\eta}}} |u_n|^{N} dx< \eta,
\end{equation}
and
\begin{equation}\label{sec}
\limsup_{n \rightarrow \infty} \int_{B^{c}_{R_{\eta}}} |u_n|^{s} dx< \eta,
\end{equation}
for all $\mu \geq M_{\eta} $. Now, as in (\ref{ple}), for $\delta>0$ and $\alpha> \alpha_0$, there exists $C>0$ such that
\begin{equation}
 |f(u_n)u_n| \leq  \delta |u_n|^N + C |u_n|^N \left(e^{\alpha|u_n|^{\frac{N}{N-1}}}-S_{N-2}(\alpha,u_n)\right).
\end{equation}

Using  H\"older's inequality and Lemma \ref{alphat11}, we get
\begin{eqnarray*}
\int_{B^{c}_{R_{\eta}}}| f(u_n)u_n | dx &\leq & \delta \int_{ B^{c}_{R_{\eta}}}|u_n|^{N}dx + C \int_{ B^{c}_{R_{\eta}}}|u_n|^N \left(e^{\alpha|u_n|^{\frac{N}{N-1}}}-S_{N-2}(\alpha,u_n)\right) dx \\
																				&\leq & \delta |u_n|^{N}_{L^N (B^{c}_{R_{\eta}})} + C |u_n|^{N}_{L^s (B^{c}_{R_{\eta}})} . 
\end{eqnarray*}
Writing
\begin{eqnarray*}
\int_{\RR^N}\left| f(u_n)u_n - f(u_\ep)u_\ep  \right|dx&\leq &  \int_{ B_{R_{\eta}}}\left| f(u_n)u_n - f(u_\ep)u_\ep  \right|dx\\ &+  &\int_{
 B^{c}_{R_{\eta}}}| f(u_n)u_n|dx+ \int_{B^{c}_{R_{\eta}}}| f(u_\ep)u_\ep  | dx 																			
\end{eqnarray*}
we conclude that 
\begin{equation}\label{ZZ1}  
\displaystyle \int_{\RR^N}\left| f(u_n)u_n - f(u_\ep)u_\ep  \right| dx= \o.
\end{equation}
Analogously
\begin{equation}\label{ZZ1}  
\displaystyle \int_{\RR^N}\left[ f(u_n)u_\ep - f(u_\ep)u_\ep  \right]dx = \o.
\end{equation}

Now, the proof follows repeating the same arguments found in \cite[Proposition 2.2]{CS}. \eproof

\subsection{Behavior of minimax levels} 

\hspace{0.6 cm} In this section, we continue studying the minimax levels. Here, we will use the notations introduced in Section 3. We also consider the Nehari manifold
\[
{\cal M}_{\mu, \ep} : = \left\{u \in E_{\mu,\ep}\setminus \{0\} \, ; \,I'_{\mu, \ep}(u)u =0\right\}
\]
and the  mountain pass minimax associated with $I_{\mu, \ep}$ given by
$$c_{\mu, \ep}= \inf \{ I_{\mu, \ep}(u) \, ;\, u \in {\cal M}_{\mu, \ep}   \}.$$

\begin{lem}\label{life} There exists $\sigma >0$, independent of $\mu$ and $\epsilon$, such that 
\begin{equation}\label{make}
||u||_{W^{1,N}(\RR^N)}^{N} \geq \sigma > 0, \,\,\, \,\,\,\mbox{for all } \, u \in {\cal M}_{\mu,\ep}.
\end{equation}
\end{lem}
\proof \, Suppose by contradiction that there exist $(\ep_n), (\mu_n) \subset (0,+\infty)$ and $(u_n) \subset {\cal M}_{\mu_n,\ep_n} $ with $||u_n|| \rightarrow 0$ as $n \rightarrow \infty$. So, there exists $ n_0 \in \mathbb{N} $ such that
\[
 ||u_n||^N < \left(\frac{\alpha_N}{\alpha_0}\right)^{N-1}, \,\,\,\,\,\mbox{for all } \, n \geq n_0.
\]
By Lemma \ref{alphat11}, there exist  $\alpha > \alpha_0$, $t_1 > 1$ and $C > 0$  independent of $n$, such that 
\begin{equation}\label{TM}
\int_{\mathbb{R}^{N}}\left(e^{\alpha|u_n|^{\frac{N}{N-1}}} -
S_{N-2}(\alpha,u_n)\right)^{t_1}dx \leq C, \,\,\,\,\mbox{for all } \, n \geq
n_0.
\end{equation}
From hypotheses $ ({H}_0) $ and $ ({H}_1) $, for each $ \eta> 0 $ and $ s > N $, there exists $ C = C ( \eta, s)> 0 $ such that
\begin{equation}
 |f(u_n)u_n| \leq  \eta |u_n|^N + C |u_n|^s \left(e^{\alpha|u_n|^{\frac{N}{N-1}}}-S_{N-2}(\alpha,u_n)\right).
\end{equation}
Using H\"{o}lder's inequality, we deduce that 
\[
||u_n||_{\mu_n,\ep_n}^{N} = \int_{\RR^N} f(u_n) u_n dx \leq  \eta |u_n|_{N}^{N} + C |u_n|^{s}_{st_2},
\]
where $t_2=\frac{t_1}{t_1-1}$.
Choosing $ \eta $ sufficiently small, it follows that
\[
||u_n||^N  \leq C_1 ||u_n||^{s}.
\]
Therefore, $||u_n||^{s-N} \geq C_2 > 0$, which is a contradiction, because $||u_n|| \rightarrow 0$.
Thus, $||u||^{N} > \sigma$, for all  $u \in {\cal M}_{\mu,\ep}$ and $\mu, \ep >0.$
 \eproof

\medskip

\begin{prop} \label{sac} For each  $ \ep \in (0,1)$,  $\displaystyle\lim_{\mu \rightarrow \infty} c_{\mu, \ep} = c_{\ep}.$
\end{prop}

\proof \, From  definitions of $c_{\mu, \ep}$ and $ c_{\ep}$, and by Remark  \ref{ab0}, 
\[
0<c_{\mu, \ep} \leq c_{\ep} < \Lambda, \,\,\,\,\, \mbox{ for all } \, \mu >0.
\]
Applying Proposition \ref{fund}, there exist sequences $\mu_n \rightarrow \infty$ and $(u_n) \subset {\cal M}_{\mu_n, \ep}$
such that
\[
I_{\mu_n, \ep}(u_n)= c_{\mu_n, \ep}   >0 \quad \mbox{and} \quad I'_{\mu_n, \ep}(u_n)=0.
\]
Then, by Proposition  \ref{rox},  there exists $u_{\ep} \in  W^{1,N}_0(\Omega_{\ep})$ solution of $(P_{\ep})$ such that, for a subsequence, we have $u_n \rightarrow u_{\ep}$ strongly in  $ W^{1,N}(\RR^N)$. Moreover, by Lemma \ref{life}, $u_{\ep} \neq 0 $ in $  W^{1,N}_0(\Omega_{\ep})$. Hence
\[
I_\ep(u_\ep)\geq c_\ep
\]
and using Proposition  \ref{rox} $iii)$
\[
c_\ep \geq \lim_{n \rightarrow \infty}c_{\mu_n, \ep} = 
\lim_{n \rightarrow \infty} I_{\ep}(u_n)= I_{\ep}(u_\ep) \geq c_\ep.
\]
Thus, 
\[
\lim_{n \rightarrow \infty}c_{\mu_n, \ep}= c_\ep.
\]
\eproof
 
As a consequence of Propositions \ref{rox} and \ref{sac}, we have the following corollary:
\begin{cor}  \label{roxe} Let $\epsilon >0$ fixed and $(u_n) \subset W^{1,N}(\RR^N)$  be a  sequence of least energy solutions of $(P_{\mu_n, \ep})$, with $\mu_{n} \rightarrow \infty$ and $\displaystyle \limsup_{n \rightarrow \infty} I_{\mu_n, \epsilon}(u_n) < \Lambda$. Then,  $(u_n)$ possesses a subsequence  that converges strongly in $ W^{1,N}(\RR^N)$ to a least  energy solution of $(P_{\ep})$.
\end{cor}

\medskip

In the sequel, let us fix $R > 2 diam (\Omega)$ such that $\Omega \subset B_R(0)$ and consider the function
$$
\xi_{\ep}(t)=
\left \{
\begin{array}{lll}
1,\quad 0\leq t \leq \frac{R}{\ep}, \\
\frac{ R}{\ep t},\quad  t \geq \frac{R}{\ep}.
\end{array}
\right.
$$
Moreover, for each $u \in W^{1,N}(\mathbb{R}^N)\setminus \{0\}$ with compact support, we set 
\[
\tilde{\beta}(u) := \frac{\int_{\mathbb{R}^N}x\xi_{\ep}(|x|)|\nabla u|^N dx}{\int_{\mathbb{R}^N}|\nabla u|^N dx}.
\]
\medskip

\begin{lem}\label {eye} There exists  $\ep^* > 0$  such that for any $\ep \in (0,\ep^*)$, there exists $\mu^*> 0$ which depends on $\ep$  such that 
if $\mu > \mu^*$, then  

$$
\tilde{\beta}(u) \in  \Omega_\e^+, \,\,\,\,\, \mbox{ for all } \,u \in {\cal M}_{\mu, \ep}  \,\,\mbox{with} \,\, I_{\mu, \ep} (u) \leq b_\ep.
$$
\end{lem}
\proof \, 
By Lemma \ref{lema 4.4}, there exists $\ep_1>0$ such that
\[
a(R,r,\ep) > c_\infty +\delta, \,\,\,\mbox{ for all } \ep < \ep_1,
\]
for some $\delta >0$.
On the other hand, by Lemma \ref{lema 4.5}, there exists $\ep_2>0$ such that
\[
b_\ep < c_\infty +\delta, \,\,\, \mbox{ for all } \ep < \ep_2.
\]
Thus, fixing $\ep^* = \min\{\ep_1,\ep_2\}$, 
\begin{equation}\label{arr}
a(R,r,\ep) > b_\ep , \,\,\, \mbox{ for all } \ep < \ep^*.
\end{equation}
Let $\ep \in (0, \ep^*) $, and suppose by contradiction that there exist sequences $\mu_n \rightarrow \infty$
and  $(u_n) \subset {\cal M}_{\mu_n, \ep}$ with $I_{\mu_n, \ep}(u_n) \leq b_{\ep}$ and
$$
\tilde{\beta} (u_n) \notin \Omega_\e^+.
$$
As $(||u_n||_{\mu_n, \epsilon})$ is bounded in $\RR$ (see Lemma \ref{lema lim}), there exists $u_{\ep} \in  W^{1,N}(\RR^N)$, such that, for a subsequence we have
\[
\left\{
\begin{array}{cl}
u_n \rightharpoonup u_{\ep}& \mbox{ in }   W^{1,N}(\RR^N),\\
u_n(x) \rightarrow  u_{\ep}(x) &\mbox{ a.e. in } \RR^N,\\
 u_n \rightarrow  u_{\ep}&  \mbox{ in }  L_{loc}^{t}(\RR^N)
\mbox{  for  }  t \geq 1.
\end{array}\right.
\]
Repeating the same arguments employed in the proof of Proposition \ref{rox}, we have that
\begin{equation}\label{ZZ1}  
\displaystyle \int_{\RR^N}\left[ f(u_n)u_n - f(u_\ep)u_\ep  \right]dx = \o.
\end{equation}
By Lemma \ref{life}, there exists $\sigma >0$ such that 
\begin{equation}\label{make}
0 < \sigma \leq ||u_n||_{W^{1,N}(\RR^N)}^{N} \leq ||u_n||_{\mu_n, \ep}^{N} = \int_{\RR^N}f(u_n)u_n dx, \,\,\, \,\,\mbox{for all } \, n \in \NN,
\end{equation}
that is, 
\begin{equation}\label{make}
0 < \sigma \leq \limsup_{n \rightarrow \infty}||u_n||_{\mu_n, \ep}^{N} = \int_{\RR^N}f(u_\ep)u_\ep dx.
\end{equation}
Hence $u_\e \neq 0$ and $I'_\e(u_\e)u_\e \leq 0$. Moreover,
$$
\lim_{n \rightarrow \infty}\tilde{\beta}(u_n) = \beta(u_\e)=y \notin \Omega_\e^{+}
$$
implying that $\Omega_\e \subset A_{\frac{R}{\e}, \frac{r}{\e},y}$. Thereby, fixing $\tau \in (0,1]$ such that $\tau u_\e \in {\cal M}_{\e , y}$, we obtain
$$
a(R,r,\e,y) \leq J_{\e,y}(\tau u_\e) \leq \liminf_{n \rightarrow \infty} I_{\mu_n,\e}(\tau u_n)  \leq \liminf_{n \rightarrow \infty} I_{\mu_n,\e}( u_n) < b_{\ep}.
$$
Using the fact that $a(R,r,\e,y)=a(R,r,\e)$,  we get a contradiction with (\ref{arr}).
\eproof

\begin{lem} \label{catego} For $0< \e <\e^*$ and $\mu > \mu^*$,
\[
cat(I_{\mu, \e}^{b_\e}) \geq cat(\Omega), 
\]
where $I_{\mu,\e}^{b_\e}:=\left\{u \in {\cal M}_{\mu,\e} \, ;\, I_{\mu, \e}(u) \leq b_\e\right\}$ and $\e^*$, $ \mu^*$ are given in Lemma \ref{eye}.
\end{lem}
\proof \,
The proof follows repeating exactly the same arguments of Lemma \ref{categ}, replacing $\beta$ by $\tilde{\beta}$ and using the Lemma \ref{eye}. 
\eproof

\medskip

\medskip

\section{Proof of Theorem 1.1 }

Let $0< \e < \e^{*}$ and $\mu > \mu^*$, where $\e^*$,$ \mu^*$ are given in Lemma \ref{eye}. Since $I_{\mu,\e}$ satisfies the $(PS)_c$ condition, for all $\mu > \mu^{*}$ and $c < b_\e$ (see Lemma \ref{lema 3.1} and Proposition \ref{fund}), we can apply the Lusternik-Schnirelman category theory and the Lemma \ref{catego} to ensure that          $I_{\mu,\e}$ has at least $cat(\Omega)$ critical points.
\eproof

\noindent{\sc Claudianor O. Alves}\\
Universidade Federal de Campina Grande,\\
Unidade Acad\^emica de Matem\'atica,\\
CEP:58429-900 - Campina Grande, PB, Brazil\\
e-mail: coalves@dme.ufcg.edu.br 

\noindent{\sc Luciana Roze de Freitas}\\
Universidade Estadual da Para\'{\i}ba,\\
Departamento de Matem\'atica,\\
CEPO:58109-790 - Campina Grande, PB, Brazil\\
e-mail: lucianarfreitas@hotmail.com


\medskip


\begin{thebibliography}{99}

\bibitem{AY} Adimurthi,  { \it Existence of positive solutions of the semilinear Dirichlet problems with critical growth for the N-Laplacian}, Ann. Sc. Norm. Super. Pisa 17 (1990), 393-413.


\bibitem{CA1} C. O. Alves,  { \it Existence and multiplicity of solutions for a class of quasilinear equation}, Adv. Nonlinear Stud. 5 (2005), 73-87.
 



\bibitem{AC6} C.O.  Alves \& G.M. Figueiredo,  { \it On multiplicity and concentration of positive solutions for a class of quasilinear problems with critical exponential growth in $\mathbb{R}^N$}, J.Differential Equations  246 (2009), 1288-1311.

\bibitem{ACL}   C.O. Alves, L.R. de Freitas \& S.H.M. Soares,  { \it Indefinite quasilinear elliptic equations in exterior domains with exponential critical growth}, Differential Integral Equations 24 (2011), 1047-1062.

\bibitem{AL}   C.O. Alves \& L.R. de Freitas,  { \it Multiplicity of nonradial solutions for a class of quasilinear equations on annulus with exponential critical growth}, Topol. Methods in Nonlinear Anal. 39 (2012), 243-262.

\bibitem{CA4} C. O. Alves  \& Y. H. Ding,   { \it Multiplicity of positive solutions to a p-Laplacian equation
involving critical nonlinearity}, J. Math. Anal. Appl. 279 (2003), 508-521.

\bibitem{CS} C. O. Alves  \& S. H. M. Soares,   { \it Multiplicity of positive solutions for a class of nonlinear Schr$\ddot{o}$dinger equations}, Adv. Differential Equation, v.15, (2010) 1083-1102.  

\bibitem{BW} T. Bartsch  \& Z. Wang ,  { \it Multiple positive solutions for a nonlinear Schr$\ddot{o}$dinger equation}, Z. Angew. Math. Phys. 51, (2000) 366-384. 

\bibitem{BC1} V. Benci  \&  G. Cerami,  { \it The effect of the domain topology on the number of positive
solutions of nonlinear elliptic problems}, Arch. Ration. Mech. Anal. 114 (1991), 79-93.


\bibitem{BC3} V. Benci  \& G. Cerami,   { \it Multiple positive solutions of some elliptic problems via the
Morse theory and the domain topology}, Cal. Var. Partial Differential Equations 02 (1994), 29-48.

\bibitem{BahC} A. Bahri  \& J.M. Coron,    { \it On a nonlinear elliptic equation involving the critical Sobolev exponent: The effect of the topology of the domain}, Commun. Pure Appl. Math. 41 (1988), 253-294.

\bibitem{JMBO}J.M. Bezerra    do \'{O}, { \it N-Laplacian equations in $\mathbb{R}^N$ with critical growth}, Abstr. Appl. Anal. 2 (1997), 301-315.

\bibitem{JMBO1}J.M. Bezerra    do \'O, { \it Semilinear Dirichlet problems for the N-Laplacian in $\mathbb{R}^{N}$ with nonlinearities in critical growth range}, Differential Integral Equations 5 (1996), 967-979.

\bibitem{OMS}J.M. Bezerra  do \'O, E.  Medeiros \&  U.  Severo,{ \it A nonhomogeneos elliptic problem involving critical growth in dimension two}, J. Math. Anal. Appl. 345 (2008), 286-304.

\bibitem{dooms} J.M. Bezerra  do \'O, E.  Medeiros \&  U. Severo, { \it On a quasilinear nonhomogeneos elliptic equation with critical growth in $\RR^N$}, J. Differential Equations 246 (2009), 1363-1386.

\bibitem{brezis} H.  Brezis, { \it An\'{a}lisis funcional. Teor\'{\i}a
y aplicaciones}, Alianza Editorial S. A., Madrid (1984).


\bibitem{Cao} D.M.  Cao, { \it Nontrivial solution of semilinear elliptic equation with critical exponent in $\mathbb{R}^2$}, Comm. Partial Differential Equations 17 (1992), 407-435.

\bibitem{BP0} G. Cerami  \& D. Passaseo,  { \it Existence and multiplicity of positive solutions for nonlinear elliptic problems in exterior domains with rich topology}, Nonlinear Anal. 18 (1992), 103-119.



















\bibitem{Ekeland}I.  Ekeland, { \it On the variational principle}, J. Math. Anal. Appl. 47 (1974), 324-353.



\bibitem{Lio} {P.L. Lions}, \textit{ The concentration-compactness principle in the calculus of
variation. The locally compact case, part 2}, Ann. Inst. H.
Poincar\'e Non Lin\'eaire, {1} (1984), 223--283.

\bibitem{moser}J.  Moser, { \it A sharp form of an inequality by N. Trudinger}, Indiana Univ. Math. J. 20 (1971), 1077-1092.


\bibitem{RP}R.  Panda, { \it On semilinear Neumann problems with critical growth for the N-Laplacian}, Nonlinear Anal 26 (1996), 1347-1366.


\bibitem{Rey} O. Rey   { \it A multiplicity result for a variational problem with lack of compactness}, Nonlinear Anal. 13 (1989), 1241-1249.



\bibitem{sergio}  E.A.B. Silva \& S.H.M. Soares, { \it Liouville-Gelfand type problems for the N-Laplacian on bounded domains of $\RR^N$}, Ann. Sc. Norm. Super. Pisa 4 (1999), 1-30.


\bibitem{tonkes}E.   Tonkes, { \it Solutions to a perturbed critical semilinear equation concerning the N-Laplacian in $\RR^N$}, Comment. Math. Univ. Carolin. 40 (1999), 679--699.

\bibitem{trudinger}N.   Trudinger, { \it On imbedding into Orlicz space and some applications}, J. Math. Mech. 17 (1967), 473-484.


\bibitem{WYZ0}Y. Wang, J.  Yang \& Y. Zhang, { \it Quasilinear elliptic equations involving the N-Laplacian with critical exponential growth in $\RR^N$}, Nonlinear Analysis 71 (2009), 6157--6169.




\end{thebibliography}
\end{document}